# CONSTRUCTING APPROXIMATIONS TO BIVARIATE PIECEWISE-SMOOTH FUNCTIONS

DAVID LEVIN

ABSTRACT. This paper demonstrates that the space of piecewise smooth functions can be well approximated by the space of functions defined by a set of simple (non-linear) operations on smooth uniform splines. The examples include bivariate functions with jump discontinuities or normal discontinuities across curves, and even across more involved geometries such as a 3-corner. The given data may be uniform or non-uniform, and noisy, and the approximation procedure involves non-linear least-squares minimization. Also included is a basic approximation theorem for functions with jump discontinuity across a smooth curve.

## 1. INTRODUCTION

High-quality approximations of piecewise-smooth functions from a discrete set of function values is a challenging problem with applications in image processing and geometric modeling. The univariate problem has been studied by several research groups, and satisfactory solutions can be found in the works of: Harten [6], Arandifa *et al.*[1], Archibald *et al.*[2, 3], Lipman *et al.*[9]. However, the 2D problem is still far from being solved, and the 1D methods are not easily adapted to the real 2D case. Furthermore, even the 1D problem is not easily solved in presence of noisy data. In the 1D problem we are given values of a piecewise smooth function, with or without noise, and the challenge is to approximate the location of the 'singular points' which separate one smooth part of the function from the other, and to also reconstruct the smooth parts. In the 2D case a piecewise smooth function on a domain $D$ is defined by a partition of the domain into segments separates by boundary curves (smooth or non-smooth), and the function is smooth in the interior of each segment. By the term smooth we mean that the derivatives (up to a certain order) of the function are bounded. Of coarse, the function and/or its derivatives may be discontinuous across a boundary curve between segments. Given data acquired from such an underlying piecewise smooth function, the challenge here is to approximate the separating curves (the singularity curves), and to reconstruct the smooth parts. Note that apart from noise in the function values, there may also be a 'noise' in the location of the separating curves (as demonstrated in Section 3.2). The problem of approximating piecewise-smooth functions is a model problem for image processing algorithms, and some sophisticated classes of wavelets and frames have been designed to approximate such functions. For example, see Candes and Donoho [5]. A method for approximation piecewise smooth functions would also be useful for the reconstruction of surfaces in CAGD or in 3D computer graphics, e.g., via the moving least-squares framework.

It is well established now that only non-linear methods may achieve optimal approximation results in 'non-smooth' spaces, e.g., see Binev *et al.*[12]. In this paper we are going back to using the 'good old' splines with uniform knots as our basis functions for the approximation, but we add to the game some (simple) non-linear operations on the space of splines. In fact, all the non-linearity used here can be expressed by the *sign* operation. We remark that the choice of the spline basis is not essential here, and other basis functions may be utilized within the same framework.

We present the idea and the proposed approximation algorithms through a series of illustrative examples. Building from derivative discontinuity in the univariate case, we move into normal discontinuity and jump discontinuity across curves in the bivariate case, with some non-trivial topologies of the singularity curves. We shall also present a basic approximation result for the case of jump discontinuity across a smooth curve. Altogether, we present a simple, yet powerful approach to piecewise smooth approximation. The suggested method seems to be quite robust to noisy data, and even the univariate version is interesting in this respect. Open issues, as the development of efficient algorithms and further approximation analysis are left for future research.

## 2. Non-smooth univariate approximations

To demonstrate the main idea, we start with the univariate problem: Assume we know that our underlying univariate function $f$ is continuous, $f \in C[a,b]$, and that it has one point of discontinuity in its first derivative in $s \in [a,b]$, and that $f'(s^-) > f'(s^+)$. Then it makes sense to look for two smooth functions $g^{[r]}$ and $g^{[\ell]}$, where $g^{[\ell]}$ approximates $f$ on the left segment $[a,s]$, and $g^{[r]}$ approximates $f$ on the right segment $[s,b]$, and such that

$$(2.1) \qquad f(x) = min(g^{[\ell]}(x), g^{[r]}(x)), \forall x \in [a,b].$$

$g^{[r]}$ may be viewed as a smooth extension of $f|_{[a,s]}$ to the whole interval $[a,b]$, and $g^{[\ell]}$ as a smooth extension of $f|_{[s,b]}$ to $[a,b]$. It is clear that there are many pairs of smooth functions $g^{[r]}$ and $g^{[\ell]}$ which satisfy the above relation. Therefore, one may suspect that the problem of finding such a pair is ill-conditioned. Let us check this by trying a specific algorithm for solving this problem, and check it on few examples. It becomes clear from these examples that the approximations by $g^{[r]}$ and $g^{[\ell]}$ are well defined in the relevant intervals, i.e., $g^{[r]}$ in $(s,b]$ and $g^{[\ell]}$ in $[a,s)$. To approximate the functions $g^{[r]}$ and $g^{[\ell]}$ we use cubic spline basis functions, with equidistant knots $t_j = a + (j-1)\delta, j = 1, ..., k, \delta = (b-a)/k$.

Assuming we are given data $\{f(x_i) | x_i \in X \subset [a,b]\}$, we look for $g^{[r]}$ and $g^{[\ell]}$ such that

$$(2.2) \qquad F_1(p) = \sum_{x_i \in X} [f(x_i) - min(g^{[\ell]}(x_i), g^{[r]}(x_i))]^2 \to minimun.$$

Here $p$ stands for the set of parameters used in the representation of the unknown functions $g^{[r]}$ and $g^{[\ell]}$. We use the convenient representation

$$(2.3) \qquad g^{[r]}(x) = \sum_{j=1}^{k} \alpha_j B_j(x), \ \ g^{[\ell]}(x) = \sum_{j=1}^{k} \beta_j B_j(x).$$

$\{B_j\}_{j=1}^{k}$ are, for example, the basis functions for cubic spline interpolant with the not-a-knot end conditions, satisfying $B_j(t_i) = \delta_{i,j}$. Hence, in (2.2) $p$ stands for the unknown splines' coefficients, $p = \{\alpha_j\}_{j=1}^{k} \cup \{\beta_j\}_{j=1}^{k}$.

In Figure 1 and in Figure 2 we see the results of reconstructing piecewise smooth functions from exact data and from noisy data. In both cases $X = \{-3 : 0.02 : 3\}$ and $\{t_j\} = \{-3 : 1.5 : 3\}$. The solution of the optimization problem (2.2) is depicted in a bold line. The underlying function $f$ is generated as $f(x) = min(f^{[\ell]}(x), f^{[r]}(x))$, and the graphs of these two generating functions are depicted by dashed lines. The fine continuous lines in the figures represent the functions $g^{[r]}$ and $g^{[\ell]}$, which, as we see in those graphs, approximate $f^{[r]}$ and $f^{[\ell]}$ accordingly, and the approximation is good only in the appropriate regions. Here $k = 5$, and thus we have 10 unknown parameters to solve for. The optimization has been performed using a differential evolution procedure, using the data values at the points $\{t_j\}$ as the starting points of the iterations for both $\{\alpha_j\}$ and $\{\beta_j\}$.

*Remark* 2.1. An alternative representation of $f$ in (2.1) is

$$(2.4) \qquad f(x) = g^{[\ell]}(x) - (g^{[\ell]}(x) - g^{[r]}(x))_+, \ \ x \in [a,b],$$

where

$$(2.5) \qquad (t)_+ = \left\{ \begin{array}{ll} t & t \geq 0 \\ 0 & t < 0 \end{array} \right. .$$

Hence, we can replace the cost functional (2.2) by

$$(2.6) \qquad F_2(p) = \sum_{x_i \in X} \left[ f(x_i) - (g_1(x_i) - (g_2(x_i))_+) \right]^2 .$$

Here $p$ stands for the set of parameters in the representation of the unknown spline functions $g_1$ and $g_2$, with the advantage that here only one unknown spline function, $g_2$, influences the functional in a non-linear manner. We shall further discuss such semi-linear cases in the bivariate case.

2.1. **The case** $f'(s^-) < f'(s^+)$ **and more.** Obviously, in this case we should replace the *min* operation within (2.2) by a *max* operation. In case we have two break points $s_1$ and $s_2$ in $[a,b]$, e.g. with $f'(s_1^-) < f'(s_1^+)$ and $f'(s_2^-) > f'(s_2^+)$, then we may look for three unknown spline functions, $g_1, g_2, g_3$, such that $min(g_1, max(g_2, g_3))$ approximates the data in the least-squares sense, and so on. To avoid high complexity we suggest to subdivide $[a,b]$ into intervals, partially overlapping, each containing at most one break point, and to blend the individual local approximations into a global one over $a, b]$. We shall further discuss and demonstrate this strategy in the 2D case.

The problem of approximating piecewise smooth univariate data has been investigated by many authors. A prominent approach to the problem is the so-called essentially non-oscillatory (ENO) and subcell resolution (SR) schemes introduced by Harten [6]. The ENO scheme constructs a piecewise-polynomial interpolant on a uniform grid which, loosely speaking, uses the smoothest consecutive data points in the vicinity of each data cell. The SR technique approximate the singularity location by intersecting two polynomials each from another side of the suspected singularity cell. In the spirit of ENO-SR many interesting works have been written using this simple but powerful idea. Recently, Arandifa *et al.*[1] gave a rigorous treatment to a variation of the technique, proving the expected approximation power on piecewise-smooth data. Archibald *et al.*[2, 3] have further improved the ENO idea by introducing polynomial annihilation techniques for locating the cell which contains the singularity. A recent paper by Lipman *et al.*[9] is using quasi-interpolation operators for this problem. Yet, the extension of the univariate methods to the 2D case is not obvious and is not simple. In [1], after locating an interval of possible singularity using ENO [6], two polynomial approximations are defined, each one approximating the data on one side of the sinularity, and their intersection is used to approximate the sinularity location. The method suggested here is similar, since we also look for two different approximations related to the two sides of a singularity. However, the least-squares optimization approach enables natural extension to interesting cases in the bivariate case. The singularity localization is integrated within the approximation procedure, and thus it is less sensitive to noise. In the next section we hope to convince that the simple idea represented in Section 2 is has the potential of solving some non-trivial bivariate approximation problems.

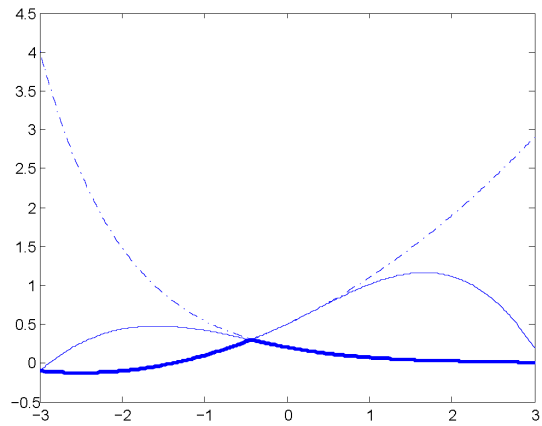

FIGURE 1.  A univariate example - No noise

## 3.  NON-SMOOTH BIVARIATE APPROXIMATIONS

As demonstrated in the 1D case, the non-linear space of functions defined by uniform splines, together with the simple operations *min* and *max*, may be used to approximated univariate piecewise smooth continuous functions. In the bivariate case we consider functions with derivative discontinuities or jump discontinuities across curves. The objectives of this section are fourfold:

(1) To exhibit a range of piecewise smooth bivariate functions which can be represented by simple non-linear operations (as *min* and *max*) on smooth functions.
(2) To suggest some non-linear least-squares approximation procedures for the approximation of piecewise smooth bivariate functions.
(3) To present interesting examples of approximating piecewise smooth bivariate functions, given noisy data.
(4) To provide a basic approximation result.

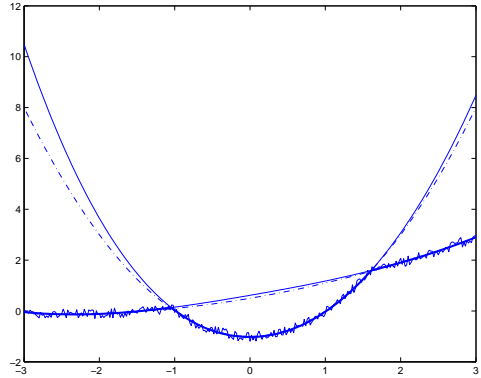

FIGURE 2.  A univariate example - Reconstruction in presence of noise

3.1. **Normals' Discontinuity across curves - Problem A.** We start with a numerical demonstration of a direct extension of the univariate approach to the approximation of continuous piecewise smooth bivariate functions. Recalling the 1D discussion, the choice of a *min* or a *max* operation depends on the sign of $f'(s^+) - f'(s^-)$. In the 2D case we refer to an analogous condition involving the slopes of the graph along the singularity curves. A discontinuity (singularity) of the normals of a bivariate function $f$ is said to be convex along a curve $\gamma$ if the exterior angle of the graph of $f$ at every point along the curve is $< \Pi$ (e.g., see Figure 6), and it is considered to be concave if the exterior angles are $> \Pi$. In a neighborhood of a concave singularity (discontinuity) curve the function may be described as the minimum between two (or more) smooth functions, and near a convex singularity curve the function may be defined as the maximum of two or more smooth functions. Let us consider the following noisy data, $\{f(x_i)\}_{x_i \in X}$, taken from a function with convex singularities. For the numerical experiment we took $X$ as the set of data points on a square grid of mesh size $h = 0.125$ in $D = [-2, 2] \times [-2, 5]$, and the given noisy data is shown in Figure 3. In this case the function has a '3-corner' type singularity, where $f$ has convex singularity along three curves meeting at a point. Therefore, we look for three spline functions, $g_1, g_2, g_3$ so that

$$(3.1) \qquad f(x) \simeq max(g_1(x), g_2(x), g_3(x)),$$

where $g_1, g_2, g_3$ solve the non-linear least-squares problem:

$$(3.2) \qquad F_A(p) = \sum_{x_i \in X} [f(x_i) - max(g_1(x_i), g_2(x_i), g_3(x_i))]^2 \rightarrow minimun.$$

Within this example we would also like to show how to blend two non-smooth approximations. Therefore, we consider the approximation problem on two partially overlapping sub-domains of $D$, $D_1 = [-2, 2] \times [-2, 2] \subset D$ and $D_2 = [-2, 2] \times [1, 5] \subset D$. After solving the approximation problem separately on each sub-domain, the two approximations will be blended into a global one. On each sub-domain the unknown functions $\{g_i\}_{i=1}^3$ are chosen to be cubic spline functions with a square grid of knots of grid size $\delta = 2$. Here again the triplet of functions $g_1, g_2, g_3$ which solve the minimization problem (3.2) is not unique. However, it turns out that the approximation to $f$ is well defined by (3.2). I.e., the parts of $\{g_i\}_{i=1}^3$ which are relevant to $max(g_1, g_2, g_3)$ are well defined.

Let us first consider the approximation on the sub-domain $D_1 = [-2, 2] \times [-2, 2] \subset D$. For the particular data shown on the left plot in Figure 4, the solution of (3.2) yields the piecewise smooth approximation depicted on the right plot. In this plot we see the full graphs of the three functions $\{g_i\}_{i=1}^3$ (for this sub-domain), while the approximation is only the upper part (the maximal values) of these graphs. The solution of the optimization problem (3.2) has been found using a differential evolution procedure [7]. As an initial guess for the three unknown functions we took, as in the univariate case, the spline function which approximates the data over the whole domain $D_1$. Next, we look for the approximation on $D_2 = [-2, 2] \times [1, 5] \subset D$, which partially overlaps $D_1$. The relevant data and the resulting approximation are shown in Figure 5.

In order to achieve an approximation over the whole domain $D = [-2, 2] \times [-2, 5]$, we now explain how to blend the two approximations defined on $D_1$ and on $D_2$. The singularity curves of the two approximations do not necessarily overlap on $D_1 \cap D_2$. Therefore, a direct blending of the two approximations will not provide a smooth transition of the singularity curve. The appropriate blending should be done between the corresponding spline functions generating these singularity curves. On each sub-domain the approximation is defined by another triplet of splines $\{g_i\}_{i=1}^3$. For the approximation over $D_2$ only two of the splines are active in the final *max* operation, and the graph of the third spline is below the maximum of the other two. To prepare for the blending step we have to match appropriate pairs of both triplets, and this can easily be done by proximity over the blending zone $D_1 \cap D_2$. The final approximation over $D$ is defined by $max(\tilde{g}_1, \tilde{g}_2, \tilde{g}_3)$, where $\{\tilde{g}_i\}_{i=1}^3$ are defined by blending the appropriate pairs, using the simplest $C^1$ blending function. The resulting blended approximation over $D$, to the data given in Figure 3, is displayed in Figure 6.

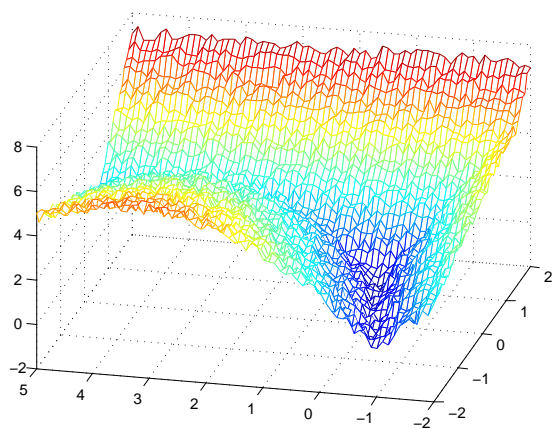

FIGURE 3. The noisy data over $[-2,2] \times [-2,5]$

3.2. **Jump Discontinuity across a curve - Problem B.** Another interesting problem in bivariate approximation is the approximation of a function with a discontinuity across a curve. Consider the case of a function defined over a domain $D$, with a discontinuity across a (simple) curve $\gamma$, separating $D$ into two sub-domains $D_+$ and $D_-$. We assume that $f|_{D_+}$ and $f|_{D_-}$ are smooth on $D_+$ and $D_-$ respectively. Such problems, and especially

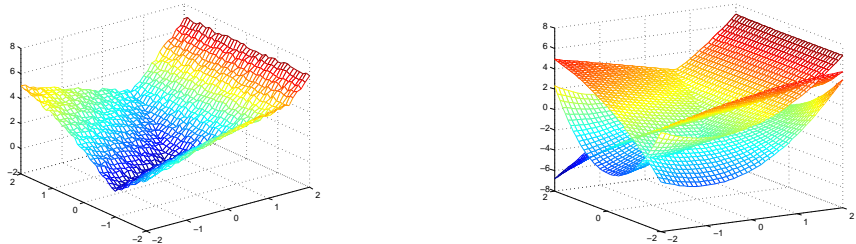

FIGURE 4. The noisy data and the 3-corner approximation over $D_1$

the problem of approximating $\gamma$, appear in image segmentation. Efficient algorithms for constructing $\gamma$, which are useful even for more involved data, are the method of snakes, or active contours, and the level-set method. The method of snakes, introduced in [8], iteratively finds contours that approach the contour $\gamma$ separating two distinctive regions in an image, with applications to shape modelling [10]. The level-set method, first suggested in [11], is also an iterative method for approximating $\gamma$, using a variational formulation for minimizing appropriate energy functionals. Recently a variational spline level-set approach has been suggested in [4]. Here, the focus is on simultaneously approximating the curve $\gamma$ and the function on $D_+$ and $D_-$. This goal is reflected in the cost functional used below, and, as demonstrated in Section 3.5, we can also handle non-simple topologies of $\gamma$, such as a 3-corner. The following procedure for treating a jump singularity comes as a natural extension of the framework for approximating a continuous function with derivative discontinuity, as suggested in Section 3.2:

Again, we look for three spline functions, $g_\gamma$, $g_+$ and $g_-$, such that the zero level set $\tilde{\gamma}$ of $g_\gamma$ approximates the singularity curve $\gamma$, $g_+$ approximates $f$ on $D_+$, and $g_-$ approximates $f$ on $D_-$. Formally, we would like to minimize the following objective function:

$$(3.3) \qquad F_B(p) = \sum_{g_\gamma(x_i)>0} [f(x_i) - g_+(x_i)]^2 + \sum_{g_\gamma(x_i)<0} [f(x_i) - g_-(x_i)]^2 \to minimun.$$

Note that the non-linearity of the minimization problem here, which we denote Problem B, is due to the non-linear operation of $sign$ checking. This approximation problem may seem to be more complicated than Problem A of the previous section, but actually it is somewhat simpler. While in problem A the unknown coefficients of all the three splines appear in a non-linear form in the objective function $F_A$ (due to the *max* operation), here only the coefficients of $g_\gamma$ influence the value of $F_B$ in a non-linear manner. This is due to the observation that once $g_\gamma$ is known, the functions $g_+$ and $g_-$ which minimize $F_B$ are defined via a linear system of equations. In view of this observation, and for reasons which will be clarified below, we use a slight variation of the optimization problem. Namely, we look for a function $g_\gamma$ which minimizes $F_B$, where $g_+$ and $g_-$ are defined by the (linear) least-squares problem:

$$(3.4) \qquad \tilde{F}_B(p) = \sum_{x_i \in X_+^h} [f(x_i) - g_+(x_i)]^2 + \sum_{x_i \in X_-^h} [f(x_i) - g_-(x_i)]^2 \rightarrow minimun,$$

where $\tilde{\gamma}$ denotes the zero level set of $g_\gamma$, $h$ is the 'mesh size' in the data set $X$, and

$$X_+^h = \{x_i \mid g_\gamma(x_i) > 0 , \; dist(x_i, \tilde{\gamma}) > h\} ,$$

$$X_-^h = \{x_i \mid g_\gamma(x_i) < 0 , \; dist(x_i, \tilde{\gamma}) > h\} .$$

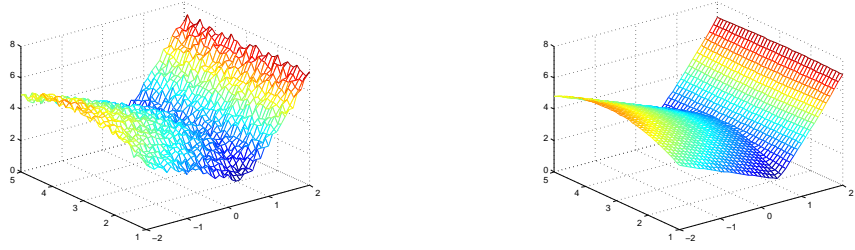

FIGURE 5. The noisy data and the approximation over $D_2$

For non-noisy data we would like to achieve an $O(h^4)$ approximation order to $f|_{D_+}$ and $f|_{D_-}$, on $D_+$ and $D_-$ respectively. This can be obtained by using proper boundary conditions in the computation of $g_+$ and $g_-$, e.g., by extending the data by local polynomial approximations. We thus consider a third version of the least-squares problem for $g_+$ and

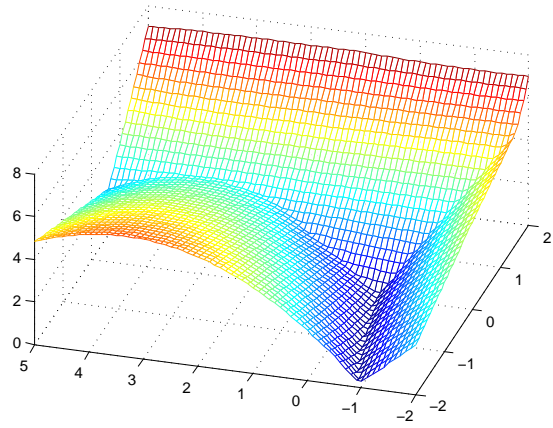

FIGURE 6. The blended approximation over $[-2,2] \times [-2,5]$

$g_-$:

(3.5)     $$\tilde{F}_B(p) = \sum_{\hat{X}_+^h} [\hat{f}_+(x_i) - g_+(x_i)]^2 + \sum_{\hat{X}_-^h} [\hat{f}_-(x_i) - g_-(x_i)]^2 \to minimun.$$

In (3.5) $\hat{X}_+^h = X \setminus X_-^h$ and $\hat{X}_-^h = X \setminus X_+^h$, $\hat{f}_+(x_i)$ is the given data $f(x_i)$ on $X_+^h$ and the extension of this data into $\hat{X}_+^h \setminus X_+^h$, and $\hat{f}_-(x_i)$ is the given data $f(x_i)$ on $X_-^h$ and the extension of this data on $\hat{X}_-^h \setminus X_-^h$. The extension operator should be exact for cubic polynomials.

*Remark* 3.1. Since $g_\gamma$ may be defined up to a multiplying factor, we may restrict its unknown coefficients to lie in a compact bounded box, and thus the existence of a global minimizer in (3.3)-(3.5) is ensured.

Let us now describe a numerical experiment based upon the above framework. The function we would like to approximate is defined on $D = [-3,3]^2$, and it has a jump discontinuity across a sinusoidal shaped curve. We may consider two types of noisy data; The first includes noise in the data values, and the second includes noise in the location of the singularity curve $\gamma$. The three unknown functions $g_\gamma, g_+, g_-$ are again cubic spline functions with a square grid of knots of grid size $\delta = 2$. However, the unknown parameters $p$ in $F_B$ are just the coefficients of $g_\gamma$. The other two spline functions are computed within the evaluation procedure of $F_B$ by solving the linear system of equations for their coefficients, i.e., the system defined by the least-squares problem (3.4). The noisy data of the second type (noise in the location of $\gamma$), and the resulting approximation obtained by minimizing (3.3), are displayed in Figures 7 and 8.

For a function with a more involved shape of singularity curve we would suggest to subdivide the domain into patches, partially overlapping, and then blend the approximations over the individual patches into a global approximation. As in the blending suggested for Problem A, the blending of two approximations to jump discontinuities over partially overlapping patches $D_1$ and $D_2$ should be performed on the functions $g_\gamma, g_+, g_-$ which generate the approximations on the different patches. Here one should take care of the fact that the function $g_\gamma$ is no uniquely defined by the optimization problem (3.3). Let us denote by $g_{\gamma,1}$ and $g_{\gamma,2}$ the functions generating the singularity curve on $D_1$ and $D_2$ respectively. In order to achieve a nice blending of the two curves we suggest to scale one of the two functions, say $g_{\gamma,1}$, so that $\alpha \cdot g_{\gamma,1} \simeq g_{\gamma,2}$ on $D_1 \cap D_2$. In fact, it is important to match the two functions only on that part of $D_1 \cap D_2$ which is close to the zero curves defined by $g_{\gamma,1}$ and $g_{\gamma,2}$.

### 3.3. Problem B - Approximation Analysis.

The approximation problem is as follows: Consider a piecewise smooth function $f$ defined over a domain $D$, with a discontinuity across a simple, smooth curve $\gamma$, separating $D$ into two open sub-domains $D_+$ and $D_-$. We assume that $f|_{D_+}$ and $f|_{D_-}$ are smooth, with bounded derivatives of order four on $D_+$ and $D_-$ respectively, and so is the curve $\gamma$. Let $X \subset D$ be a grid of data points of grid size $h$, and let us consider the approximations for Problem B using bi-cubic spline functions with knots on a grid of size $\delta = mh$ $(m > 3)$. The classical result on approximation by least-squares by cubic splines implies an $O(h^4)$ approximation order to a function with bounded derivatives of order four (provided there are enough data points for a well-posed solution). On the other hand, even in the univariate case, the location of a jump discontinuity in a piecewise smooth function is inherently up to an $O(h)$ error. Therefore, the best we can expect from a good approximation procedure for $f$ such as above is the following:

**Theorem 3.2.** *Consider Problem B on $D$ and let $g_\gamma$ be a bi-cubic spline function (with knots' grid size $\delta = mh$) which gives a local minimum to (3.3), with $g_+$ and $g_-$ defined by minimizing (3.5). Denote the segmentation defined by $g_\gamma$ by $G_+ = \{\, g_\gamma(x_i) > 0 \,,\, x_i \in X \,\}$ and $G_- = \{\, x_i \mid g_\gamma(x_i) < 0 \,,\, x_i \in X \,\}$. For $C > 0$, and for $h$ small enough, there exists such local minimizer $g_\gamma$ such that if $x_i \in G_+ \cap D_-$ or $x_i \in G_- \cap D_+$ then $dist(x_i, \gamma) < Ch^3$.*

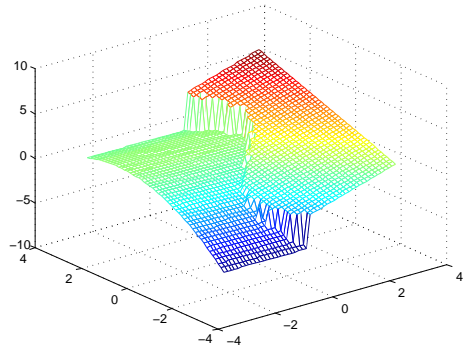

FIGURE 7. Discontinuity across a noisy curve $[-3,3] \times [-3,3]$

*Proof.* The Theorem says that the zero level set of $g_\gamma$, $\tilde{\gamma}$, separates well the data set $X$ into the two parts, and only data points which are very close to $\gamma$ may appear in the wrong segment. In order to prove this result we first observe that the curve $\gamma$ can be approximated by the zero level set of bi-cubic splines with approximation error $< C_1 h^4$. One such spline would be $s_\gamma$, the approximation to the signed distance function related to the curve $\gamma$. Fixing $g_\gamma = s_\gamma$ determines $g_+$ and $g_-$ which minimize $\tilde{F}_B$ for this $g_\gamma$, and we denote the corresponding value $F_B[s_\gamma]$. We note that the contribution to the value of $F_B$ is $o(h^4)$ (as $h \to 0$) from a point which fall on the right side of $\tilde{\gamma}$, and it is $O(1)$ from a point on the wrong side of $\tilde{\gamma}$. For a small enough $h$, only a small number of points $x_i \in X$ will fall in the wrong side of $\tilde{\gamma}$, and any choice of $g_\gamma$ which induces more points in the wrong side will induce a larger value of $F_B$. Obviously, the minimizing solution induces a value $F_B \leq F_B[s_\gamma]$ , and this can be achieved only by reducing the set of 'wrong side' points. Since $g_\gamma = s_\gamma$ already defines an $O(h^4)$ separation approximation, only points which are at distance $O(h^4)$ from $\gamma$ may stay on the wrong side in the local minimizer which evolves by a continuous change of $s_\gamma$ which reduces $F_B$.    □

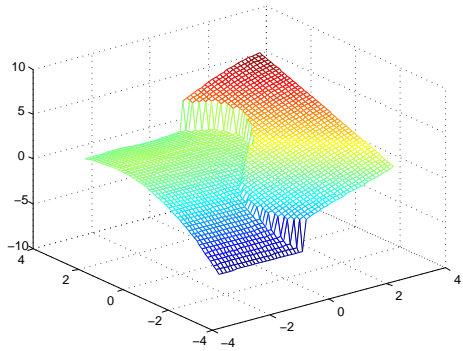

FIGURE 8. The approximation using noisy curve data

**Corollary 3.3.** *If the least-squares problems defining $g_+$ and $g_-$ by (3.4) are well-posed, we get*

$$(3.6) \qquad \|f - g_+\|_{\infty, D_+} \leq C_2 h^4,$$

$$(3.7) \qquad \|f - g_-\|_{\infty, D_-} \leq C_3 h^4.$$

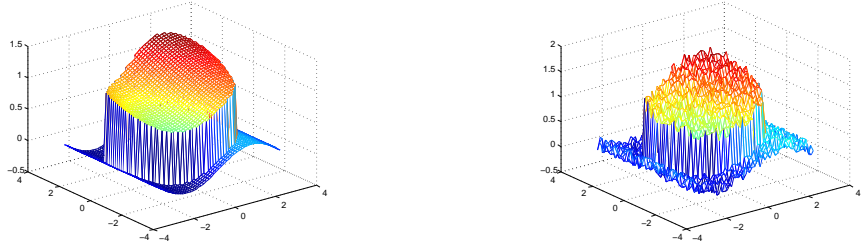

FIGURE 9. The underlying function and its 1st type noisy data.

*Remark* 3.4. The above well-posed condition can be checked while computing $g_+$ and $g_-$. Also, an $O(h^4)$ approximation order can be obtained by using proper boundary conditions in the computation of $g_+$ and $g_-$, e.g., by extending the data by local polynomial approximations, as suggested in (3.5).

*Remark* 3.5. The need to restrict the set of data points defining $g_+$ and $g_-$ in (3.4) emerged in view of the condition needed for the proof of Theorem 3.2. As shown in the numerical example below, this restriction may be very important in practical applications.

3.4. **Noisy data of the 1st type.** This section demonstrate the performance of the method for the approximation of noisy data of a function with jump discontinuity. Furthermore, we use this example to emphasize the importance of using the restricted sets in $\tilde{F}_B$ rather than using $F_B$. The underlying function and its noisy version are displayed in Figure 9. In the numerical test we have used the same mesh and knot sizes as in the previous example. In figure 10 we show the results with and without restricting the the set of points which participate in the computation of $F_B$. In the left graph we note that the approximation in the inner region is infected by wrong values from the outer region, and this is clearly corrected in the right graph where the least-squares approximations use values which are not to close the discontinuity curve. In Figure 11 we see two approximations to the exact singularity curves (in red), using different knots' grid sizes, $\delta = 1.5$ and $\delta = 2.$, together with the singularity curve of the underlying function. As expected, a smaller $\delta$ enables higher flexibility of $g_\gamma$ and a better approximation to the exact curve.

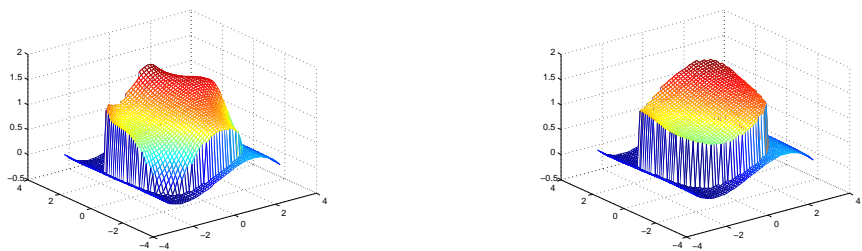

FIGURE 10. Two approximations using different point sets in the least-squares approximation.

3.5. **3-corner Jump Discontinuity - Problem C.** Combining elements from Problems A and B, we can now approach the more complex problem of a 3-corner discontinuity. Consider the case of a function defined over a domain $D$, with a discontinuity across three curves meeting at a 3-corner, subdividing $D$ into three sub-domains $D_1$, $D_2$ and $D_3$, as in Figure 12. We assume that $f|_{D_i}$ is smooth on $D_i$, i=1,2,3. Following the above discussions, the following procedure is suggested:

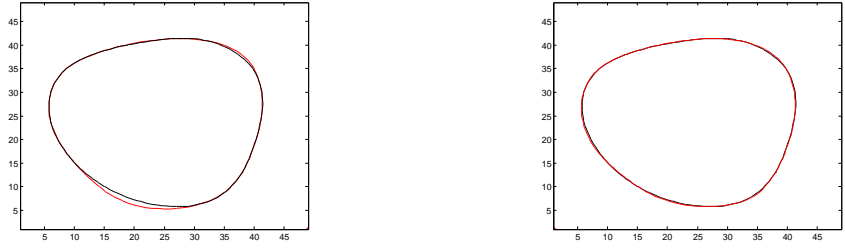

FIGURE 11. Two approximations to the exact singularity curves - using different knots' grid sizes, $\delta = 2$. and $\delta = 1.5$.

We look for three spline functions, $\{g_i\}_{i=1}^3$, approximating $f$ on $\{D_i\}_{i=1}^3$ respectively. Here, the approximation of the segmentation into three domains cannot be done via a zero level set approach. Instead, we look for an additional triplet of spline functions, $\{h_i\}_{i=1}^3$ which define approximations $\{E_i\}_{i=1}^3$ to $\{D_i\}_{i=1}^3$ as follows:

$$E_1 = \{ x \mid h_1(x) > max(h_2(x), h_3(x) \},$$
$$E_2 = \{ x \mid h_2(x) > max(h_1(x), h_3(x) \},$$
$$E_3 = \{ x \mid h_3(x) > max(h_1(x), h_2(x) \}.$$

Denoting $\|u\|_{2,E} = \sum_{x_i \in E} u^2(x_i)$ we would like to minimize the following objective function:

$$(3.8) \qquad F_C(p) = \|f - g_1\|_{2,E_1} + \|f - g_2\|_{2,E_2} + \|f - g_3\|_{2,E_3} \to minimun.$$

Hence, the segmentation is defined by a *max* operation as in Problem A. Given a segmentation of $D$ into $\{E_i\}_{i=1}^3$, the triplet $\{g_i\}_{i=1}^3$ is defined, as in Problem B, by a system of linear equations which defines the least-squares solution of (3.8). To achieve better approximation on $\{D_i\}_{i=1}^3$, in view of Theorem 3.2, the least-squares approximation for $\{g_i\}_{i=1}^3$ should exclude data points which are near joint boundaries of $\{E_i\}_{i=1}^3$.

For a numerical illustration of Problem C and the approximation obtained by minimization of $F_C$ we took noisy data from a function with 3-corner discontinuity in $D = [-2,2]^2$. All the unknown spline functions, $\{g_i\}_{i=1}^3$ and $\{h_i\}_{i=1}^3$ are bi-cubic with a square grid of knots of grid size $\delta = 2$. Since only the splines $\{h_i\}_{i=1}^3$ enter in a non-linear way into $F_C$, the minimization problem involves $3 \times 9 = 27$ unknowns. As in all the previous examples we have used a differential evolution algorithm for finding an approximate solution of this minimization problem. The noisy data and the resulting approximation are shown in Figure 12.

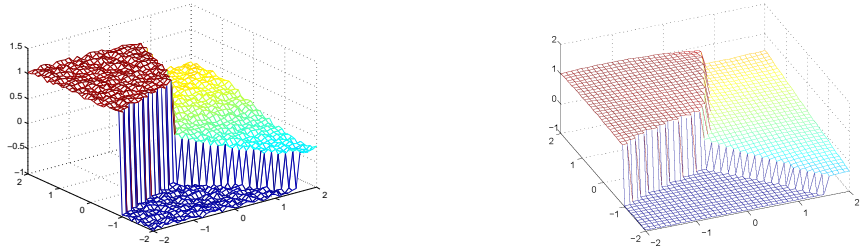

FIGURE 12. 3-corner discontinuity - noisy data and approximation.

## 4. SUMMARY AND ISSUES FOR FURTHER RESEARCH

We have introduced a unified framework for approximating functions with normals' discontinuity or jump discontinuity across curves. The method may be viewed as an extension of the well known procedures of boolean operations in solid geometry. In this work it is suggested to use a kind of boolean operations on splines as an approximation tool. Through a series of non-trivial examples we have presented the potential of this approach to achieve high quality approximations in many applications. It is interesting to note that all the non-linearity in the suggested approximations can be expressed by the *sign* operation, or equivalently the $(\cdot)_+$ operation. The approximation procedure requires high dimensional non-linear optimization, and thus the complexity of computing the approximations is very high. For all the numerical examples in this paper we have used a very powerful matlab code written by Markus Buehren, based upon the method of Differential Evolution ([7]). The execution time ranges from 1 second for the simple univariate problem to 80 seconds for the bivariate Problem C. The differential evolution algorithm usually finds a local minima and not 'the global minimizer'. Yet, as demonstrated, is finds for us very good approximations, and it seems to be robust to noise. A main issue for further study would be the acceleration of the optimization process, e.g., by generating good initial candidates for the optimization. Yet, in spite of the high computational cost, the method may still be very useful for high quality up-sampling, and for functions (or surfaces) with few singularity curves. In a scene with many discontinuities we would suggest to subdivide the domain into patches, each containing one or two singularity curves. Choosing partially overlapping patches, the local approximations can be blended into a global approximation, as demonstrated in Section 3.2. Another simple idea is based upon Corollary 3.3, which tells us to ignore few data points near the approximated singularity curve, in order to attain higher approximation order.

Other important issues for further research would be the following:

(1) Improved optimization: Here we believe that geometric considerations may be used to significantly accelerate the minimization procedure. Gradient-descent algorithms, similar to those used in [4], may also be helpful here.
(2) Simple rules for choosing grid size for the splines.
(3) Other basis functions instead of splines.
(4) Using $\ell_1$-norm instead of the $\ell_2$-norm in the objective functions.
(5) Application to 3D surface data via the moving least-squares method.
(6) Further approximation analysis.